\documentclass{article}
\usepackage[cp1252]{inputenc}
\usepackage{amsmath,amsfonts,amscd,amssymb,amsthm,mathrsfs}
\usepackage[mathcal]{eucal}
\newtheorem{lemma}{Lemma}[section]
\newtheorem{prop}[lemma]{Proposition}
\newtheorem{theorem}[lemma]{Theorem}

\newtheorem*{Rq}{Remark}
\newtheorem{Cor}{Corollary}
\def\Hom{\mathop{\rm Hom}\nolimits}

\newcommand{\Gal}{\mathop{\rm Gal}\nolimits}

\newcommand{\zp}{\mathbb{Z}_p}
\newcommand{\qp}{\mathbb{Q}_p}
\newcommand{\on}{\mathcal{O}_n}
\def\Frac{\mathop{\rm Frac}\nolimits}
\def\|{\mathop{\rm |}\nolimits}
\newcommand{\limproj}{\mathop{ {\lim\limits_{\longleftarrow}}} \nolimits}

\def\entL{{[\mkern-2.5mu[}}
\def\entR{{]\mkern-2.5mu]}}
\author{\textsc{Bruno Angl\`es, Thomas Herreng}}
\title{\textbf{On a result of Iwasawa}}
\date{}

\begin{document}
\maketitle
\begin{center}
Universit\'e de Caen
\\Laboratoire de Math\'ematiques Nicolas Oresme
\\CNRS UMR 6139
\\Boulevard Mar\'echal Juin
\\B.P. 5186
\\14032 Caen Cedex
\\FRANCE
\end{center}
\bigskip
E-mail: bruno.angles@math.unicaen.fr, thomas.herreng@math.unicaen.fr

\abstract{ We recover a result of Iwasawa on the $p$-adic logarithm of principal units of $\qp(\zeta_{p^{n+1}}$
by studying the value at $s=1$ of $p$-adic $L$-functions.}

\\Math Subject classification: 11R23, 11R18.

\bigskip
Since the ninteenth century, it is known that values of $L$-functions at $s=1$ contain 
deep arithmetic information. This result has much importance because it links analytic formulas
with arithmetic invariants. Kubota and Leopoldt have defined an analogous $L$-function on the $p$-adic
fields with analytic techniques. Iwasawa has shown how to construct this $p$-adic $L$-function algebraically.
Our aim in this paper is to give some algebraic interpretations to the analytic formulas giving the 
value at $s=1$ of $p$-adic $L$-functions. It leads us to study some properties of the $p$-adic
logarithm which enable us to establish the Galois module structure of the logarithm of principal units.
The result obtained (theorem 1.10 below) had been discovered first by Iwasawa (\cite{IWA3}) in 1968 
using explicit reciprocity laws. The first corollary we state to this theorem is an important 
result also due to Iwasawa (\cite{IWA2}) which gives the structure of the plus-part of the principal 
units modulo cyclotomic units.

In the second section we use the theorem 1.10 
to study the minus part of the projective limit of these units. The aim is to obtain a theorem
which looks like the one of Iwasawa (\cite{IWA2}).
But in the minus-part there are no more cyclotomic units. Yet by considering the $p$-adic 
logarithm of principal units we are able to obtain an Iwasawa-like result (theorem 2.3) and its
global counterpart.

The authors would like to thank Maja Volkov for her careful reading of the paper and great 
improvements in the presentation.

For all the paper, fix an odd prime number $p$ and denote by $v_p$ the normalized valuation
at $p$. We are interested in cyclotomic $p$-extensions. We start with $\theta$ an even character
of $\Gal(\mathbb Q(\mu_{pd})/\mathbb Q)$ of conductor $f_\theta=d$ or $pd$ with $d \geq 1$ and $p \nmid
d$. For $n \geq 0$, let $q_n=p^{n+1}d$. For all integer $m$, we identify both groups
$$\begin{array}{lcl}\Gal(\mathbb Q(\mu_m)/\mathbb Q) &\simeq& (\mathbb Z/m\mathbb Z)^\times\\
\sigma_a &\leftrightarrow & a \quad \quad \quad \quad \mbox{ where $\sigma_a(\zeta_m)=\zeta_m^a$}.
\end{array}$$
Fix a primitive $d$th root
of unity $\alpha$. We note $\zeta_{q_n}=\alpha \zeta_{p^{n+1}}$
with $\zeta_{p^{n+1}}$ defined such that
$\zeta_{p^{n+1}}^p=\zeta_{p^n}$ and $\zeta_p \in \mu_p \setminus
\{1\}$. Let $F$ be the Frobenius endomorphism of $\qp(\alpha)/\qp$
(i.e. $F(\alpha)=\alpha^p$), and $\Delta=\Gal(\mathbb
Q(\mu_{pd})/\mathbb Q)$. Then $F$ corresponds to $\sigma_p \in
\Delta$. As $\theta \in \widehat{\Delta}$ we write
$\theta(F)$ for $\theta(\sigma_p)$.\\ \indent Let
$K_n=\qp(\mu_{q_n})$ and $K_\infty = \bigcup_{n\geq0} K_n$. We
know that $\Gamma_n=\Gal(K_n/K_0) \simeq \mathbb Z/p^n \mathbb Z$
and $\Gamma=\Gal(K_\infty/K_0) \simeq \zp$. Moreover
$\gamma_0=\sigma_{1+q_0}$ is a topological generator of $\Gamma$.
We have the decomposition
$$\begin{array}{lcl}\Gal(K_n/\mathbb Q) &\simeq& \Delta \times \Gamma \\
\sigma_a &\mapsto& (\delta(a), \gamma_n(a)).
\end{array}$$

\indent Finally, we fix some notations from Iwasawa theory. Let
$\mathcal O$ be the ring of integers of $\qp(\theta)$ and
$\Lambda=\mathcal O\entL T\entR$. By a fundamental results of
Iwasawa there exists a power series $f(T, \theta)\in
\Lambda$ such that for all $n \geq 0$, $\chi \in
\widehat{\Gamma}_n$, and $s \in \zp$ $$L_p(s, \theta \chi)=
f(\overline{\chi}(\gamma_0) (1+q_0)^s-1, \theta).$$

\section{Value at $s=1$ of $p$-adic $L$-functions}

\subsection{A Lemma on the Iwasawa algebra}
Let $\Lambda=\zp\entL T\entR$ be the Iwasawa algebra, $g \in
\Lambda$ a formal power series and $\omega_n$ the element
$(1+T)^{p^n}-1$. It is known that the Iwasawa algebra can be
described as follows. Let $\Gamma$ be a multiplicative
topological group isomorphic to $\zp$, $\Gamma_n$ its quotient
$\Gamma/\Gamma^{p^n}\simeq \mathbb Z/ p^n \mathbb Z$, and
$\zp\entL\Gamma\entR =\limproj \zp[\Gamma_n]$. Let $\gamma_0$ be a topological generator of
$\Gamma$ and $\overline{\gamma_0}$ its image in $\zp[\Gamma_n]$.We thus have an isomorphism
$$\begin{array}{ccl} \zp[\Gamma_n] & \simeq &
\zp[T]/((1+T)^{p^n}-1) \\ \overline \gamma_0 & \mapsto &\ 1+T \end{array}$$
which, passing to the projective limit, gives an isomorphism $$\begin{array} {ccl}
\zp\entL\Gamma\entR & \simeq & \zp\entL T\entR\\ \gamma_0 & \mapsto&
\ 1+T \end{array}$$ (see for example \cite{WAS}). Via this
isomorphism the power series $g$ can be written as a sequence of
elements $\epsilon_n \in \zp[\Gamma_n]$ compatible with the
restriction morphisms. The aim of this paragraph is to express
$\epsilon_n$ according to $g$. \\ \indent We can write $g$ as $$
g(T)= \sum_{b=0}^{p^n-1} a_n(b) (1+T)^b + \omega_n(T)Q_n(T)$$ with
$Q_n(T) \in \Lambda$. 
The canonical isomorphism described above implies that $\epsilon_n =
\sum_{b=0}^{p^n-1} a_n(b) \overline{\gamma_0}^b$. Pick a
character $\chi \in \widehat{\Gamma}_n=\Hom(\Gamma_n, \mathbb Q_p^\times)$, and let $e_\chi = \frac 1{p^n} \sum_{\gamma \in \Gamma_n}
\overline{\chi}(\gamma) \gamma\in
\qp[\Gamma_n]$ be the associate idempotent. Thus

 \begin{eqnarray*}
e_\chi \epsilon_n & = & e_\chi ( \sum_{b=0}^{p^n-1} a_n(b) \overline{\gamma_0}^b) \\
 &=& (\sum_{b=0}^{p^n-1} a_n(b) \chi(\gamma_0)^b) .e_\chi\\
  &=& g(\chi(\gamma_0)-1)e_\chi .
\end{eqnarray*}
Summing over all characters of $\Gamma$
we deduce the following result.
\begin{lemma} Let $g$ be a power series in $\zp\entL T\entR$ and $\epsilon_n$
its image in $\zp[T]/\omega_n(T) \zp[T]$. Via the canonical isomorphism
between $\zp\entL T\entR$ and $\zp\entL\Gamma\entR$ the element $\epsilon_n$
in $\zp[\Gamma_n]$ writes as
$$\epsilon_n = \sum_{\gamma \in \Gamma_n} \left[ \frac{1}{p^n} \sum_{\chi \in \widehat{\Gamma}_n}
g(\chi(\gamma_0)-1) \overline{\chi}(\gamma) \right] \gamma$$ where $\gamma_0$
is a topological generator of $\Gamma$.
\end{lemma}

\subsection {Value of the $p$-adic $L$-function at $s=1$ : algebraic interpretation}
Our aim is to give an algebraic interpretation of the following formula
(which can be found in \cite{WAS}, theorem 5.18) :
\begin{prop} Let $\chi$ be an even nontrivial character of conductor $f$ and $\zeta$ be a primitive
$f$th root of unity. We have
$$L_p(1, \chi)=-(1-\frac{\chi(p)}{p}) \frac{\tau(\chi)}{f} \sum_{a=1}^f \overline{\chi}(a) \log_p(1-\zeta^a)$$
where $\tau(\chi)$ is the Gauss sum $\tau(\chi)=\sum_{a=1}^f \chi(a) \zeta^a$.
\end{prop}

We use the notations of the introduction.
Recall the well-known relation for all $n \geq 0$ and $i \leq n$
\begin{equation} N_{K_n/K_i}(\zeta_{q_n}^{F^{i-n}}-1) = \zeta_{q_i}-1 .
\end{equation}

Assume $\theta \neq 1$. By the results of the previous paragraph we know there exists
$(\epsilon_n(\theta))_{n \geq0} \in \limproj \mathcal O[\Gamma_n]$
corresponding to $f(\frac{1+q_0}{1+T}-1, \theta)$ and that
for all $n \geq 0$ and $\chi \in \widehat{\Gamma}_n$
$$e_\chi \epsilon_n = L_p(1, \theta \chi) e_\chi . $$
\begin{Rq}Assume $\theta \neq 1$. We have for all $\chi \in \widehat{\Gamma}_n$ and $n \geq 0$
$L_p(1, \theta \chi) \neq 0$. Thus for all $n\geq 0$ $\epsilon_n \in
\Frac(\mathcal O)[\Gamma_n]^\times$.
\end{Rq}

Let $\chi \in \widehat{\Gamma}_n$. We assume that $\chi \neq 1$ if $f_\theta=d$. Then the
conductor of the product $\theta \chi$ is $f_{\theta \chi}=q_k$ for an integer $0 \leq k \leq n$.
\begin{lemma}We have
$$\frac{1}{p^n} \theta(F^{n-k}) \tau(\overline{\theta \chi}) L_p(1, \theta \chi)=
- \sum_{\delta \in \Delta} \overline \theta (\delta) e_\chi \log_p(1-\zeta_{q_n}^\delta).$$
\end{lemma}

\begin{proof}
Start from the formula in proposition 1.2. Using the fact that $\tau(\chi) \overline{\tau(\chi)} =f$
for a character $\chi$
of conductor $f$ we find
$$\tau(\overline{\theta \chi}) L_p(1, \theta \chi)=
-\sum_{\delta \in \Delta, \gamma \in \Gamma_k} \overline \theta (\delta) \overline{\chi}(\gamma) \log_p(1-\zeta_{q_k}^{\delta \gamma}).$$
Then using the equality (1)
\begin{eqnarray*}\tau(\overline{\theta \chi}) L_p(1, \theta \chi) & = &
-\sum_{\delta \in \Delta, \gamma \in \Gamma_n} \overline \theta (\delta) \overline{\chi}(\gamma) \log_p(1-\zeta_{q_k}^{\delta \gamma F^{k-n}})\\
&=&- \sum_{\delta \in \Delta, \gamma \in \Gamma_n} \overline \theta (\delta)\overline{\theta}(F^{n-k}) \overline{\chi}(\gamma) \log_p(1-\zeta_{q_n}^{\delta \gamma})\\
&=&-p^n \overline{\theta}(F^{n-k}) \sum_{\delta \in \Delta} \overline \theta (\delta) e_\chi \log_p(1-\zeta_{q_n}^{\delta })
\end{eqnarray*}
which proves the lemma.
\end{proof}
Assume now that $f_{\theta}=pd$ and let $\dot{\mathscr T}_n= \sum_{i=0}^n F^{n-i} p^{i-n} \zeta_{q_i}$.
\begin{prop} For all $n \geq 0$ we have
$$ \sum_{\delta \in \Delta} \overline{\theta}(\delta) \epsilon_n(\theta) (\dot{\mathscr T}_n^\delta)=
-\sum_{\delta \in \Delta} \overline{\theta}(\delta) \log_p(1-\zeta_{q_n}^\delta).$$
\end{prop}
\begin{proof}
Recall that for a character $\chi$ of second kind (i.e. whose conductor is $f_\chi = p^{k+1}$)
and an integer $l$ $0\leq l\leq n$
the $\chi$-part $e_\chi \zeta_{p^{l+1}}$ vanishes if and only if $l \neq k$ (see for example \cite{Erez}).
Let $\chi$ be a character of $\Gamma_n$. We have already shown that
$\sum_{\delta \in \Delta} \overline{\theta}(\delta) e_\chi \epsilon_n(\theta)(\dot{\mathscr T}_n^\delta)
= \sum_{\delta \in \Delta} \overline{\theta}(\delta) L_p(1, \theta\chi) e_\chi (\dot{\mathscr T}_n^\delta)$.
Moreover
$$e_\chi(\dot{\mathscr T}_n^\delta)=\frac{1}{p^n} \sum_{\gamma \in \Gamma_n}
p^{k-n}\overline{\chi}(\gamma) \zeta_{q_k}^{\delta \gamma F^{n-k}}$$
where $k$ is such that $f_\chi=p^{k+1}$ when $\chi \neq 1$ and $k=0$ when $\chi =1$.
Thus we have  $0 \leq k \leq n$ and $f_{\theta \chi}=q_k$. Then
\begin{eqnarray*}
\sum_{\delta \in \Delta}e_\chi \overline{\theta}(\delta) \epsilon_n(\theta) (\dot{\mathscr T}_n^\delta) &=&
\frac{L_p(1, \theta\chi)}{p^n} \sum_{\stackrel{\delta\in \Delta}{\gamma \in \Gamma_n}}
p^{k-n}\overline {\theta}(\delta) \overline{\chi}(\gamma)\zeta_{q_k}^{\delta\gamma F^{n-k}} \\
&=& \frac{\theta(F^{n-k}) L_p(1, \theta \chi)}{p^n}
\sum_{\stackrel{\delta\in \Delta}{\gamma \in \Gamma_n}}
p^{k-n}\overline {\theta}(\delta) \overline{\chi}(\gamma)\zeta_{q_k}^{\delta\gamma} \\
&=& \frac{\theta(F^{n-k}) L_p(1, \theta \chi)}{p^n}
\sum_{\stackrel{\delta\in \Delta}{\gamma \in \Gamma_k}}
\overline {\theta}(\delta) \overline{\chi}(\gamma)\zeta_{q_k}^{\delta\gamma} \\
&=&\frac{1}{p^n} \theta (F^{n-k}) L_p(1, \theta \chi) \tau(\overline{\theta\chi})\\
&=& -\sum_{\delta \in \Delta} \overline{\theta}(\delta)e_\chi \log_p(1-\zeta_{q_n}^\delta).
\end{eqnarray*}
The last equality follows from the previous lemma. Summing over all characters $\chi \in
\widehat{\Gamma}_n$ we obtain the required equality.
\end{proof}

Of course there is a similar equality in the case where $f_\theta = d$ instead of $pd$. However
the result is slightly more complicated to state. We need the following lemma.
\begin{lemma}If $f_\theta =d$ then
\begin{itemize}
\item $\sum_{\delta \in \Delta} \overline{\theta}(\delta) \log_p(1-\zeta_{q_0}^\delta)=
(\theta(F)-1) \sum_{y \in (\mathbb Z/d \mathbb Z)^\times} \overline{\theta}(y) \log_p(1- \alpha^y) \mbox{\rm{ and}}$
\item $\tau(\overline{\theta})= -\sum_{\delta \in \Delta} \overline{\theta}(\delta) \zeta_{q_0}^\delta$
\end{itemize}
\end{lemma}
It follows that for all n $\geq 0$,
$$\frac{1}{p^n}\tau(\overline{\theta}) L_p(1,\theta) =
-(1-\frac{\theta(F)}{p}) \sum_{\delta \in \Delta} \overline{\theta}(\delta) e_{\chi_0}
\log_p(1-\zeta_{q_n}^\delta)$$ where $\chi_0$ is the trivial character of $\Gamma_n$. We can now state the
theorem which includes all the results of this section.

\begin{theorem}Let $\theta$ be an even non-trivial character of $\Gal(\mathbb Q(\mu_{pd})/\mathbb Q)$
and $(\epsilon_n(\theta))_{n\geq 0} \in \limproj \mathcal
O[\Gamma_n]$ the element corresponding to $f(\frac{1+q_0}{1+T}-1,
\theta)$. We define $$\dot{\mathscr T}_n = \left\{
\begin{array}{ll} \sum_{i=0}^n F^{n-i} p^{i-n} \zeta_{q_i} &
\textrm{ when $f_\theta=pd$}\\ \sum_{i=1}^n F^{n-i} p^{i-n}
(1-\frac F p )\zeta_{q_i}-p^{-n}F^n(f-1) \zeta_{q_0}& \textrm{
when $f_\theta=d$}.\end{array} \right.$$ We thus have for all
$n\geq 0$, $$\sum_{\delta \in \Delta} \overline{\theta}(\delta)
\epsilon_n(\theta) (\dot{\mathscr T}_n^\delta)=
-E(\theta)\sum_{\delta \in \Delta} \overline{\theta}(\delta)
\log_p(1-\zeta_{q_n}^\delta)$$ where $E(\theta)$ is a kind of
Euler factor : $E(\theta)=(1-\frac{\theta(F)}p)$ when $f_\theta=d$
and $E(\theta)=1$ when $f_\theta=pd$.
\end{theorem}

\subsection{An application}

We now apply theorem 1 to the case where $f_\theta=pd$. We
slightly change our notations. For the rest of the paper we denote
by $\Delta$ the group $(\mathbb Z/p \mathbb Z)^\times$, $\theta_1$
a character of $\Delta$ and $\theta_2$ a character of conductor
$d$ with $d$ dividing $p-1$ such that $\theta=\theta_1 \theta_2$
is even. We assume first that $\theta_1 \neq 1, \omega$. Both
cases will be treated separately in the sequel. Let
$W=\zp[\theta_2]$ and $(\epsilon_n(\theta)) _{n \geq 0} \in W\entL
T \entR$ the element corresponding to $f(\frac{1+pd}{1+T}-1,
\theta)$. Recall that, as usual, $\alpha$ is a primitive $d$th
root of unity. In this situation we have the following result.
\begin{lemma}
Let $\mathscr T_n=\sum_{i=0}^n p^{i-n} \zeta_{p^{i+1}}$. Then
$$\tau(\overline \theta_2) \epsilon_n(\theta) e_{\theta_1} \mathscr T_n =
-\sum_{y\in (\mathbb Z/d \mathbb Z)^\times} \overline{\theta }_2(y) e_{\theta_1} \log_p(\alpha^y - \zeta_{p^{n+1}}).$$
\end{lemma}
\begin{proof}
Theorem 1 can be restated as $$\sum_{\delta \in (\mathbb Z/pd
\mathbb Z)^\times} \overline \theta (\delta)
\epsilon_n(\theta)(\dot{\mathscr T}_n^\delta)=-\sum_{\delta \in
(\mathbb Z/pd \mathbb Z)^\times} \overline \theta(\delta)
\log_p(1-\zeta_{q_n}^\delta)$$ where $\dot{\mathscr T}_n=
\sum_{i=0}^n p^{i-n} \zeta_{q_i}$ and $q_i=dp^{i+1}$. A
straightforward calculation gives
\begin{eqnarray*} \sum_{\delta \in (\mathbb Z/pd \mathbb Z)^\times} \overline \theta (\delta)
\epsilon_n(\theta)(\dot{\mathscr T}_n^\delta) &=&
\sum_{\stackrel{\delta \in (\mathbb Z/p \mathbb Z)^\times}{y \in (\mathbb Z/d \mathbb Z)^\times}}
\overline \theta_2 (y) \overline \theta_1(\delta)\epsilon_n(\theta)(\sum_{i=0}^n p^{i-n} \alpha ^y \zeta_{p^{d+1}}^\delta)\\
&=& (p-1) \sum_{y \in (\mathbb Z/d \mathbb Z)^\times} \overline \theta_2 (y) \alpha^y
\epsilon_n(\theta) e_{\theta_1}(\sum_{i=0}^n p^{i-n} \zeta_{p^{i+1}})\\
&=& (p-1) \tau(\overline \theta_2) \epsilon_n e_{\theta_1} \mathscr T_n.
\end{eqnarray*}
\end{proof}
We need a relation between primitive characters and unprimitive
ones which is given by the following lemma.
\begin{lemma}Let $\theta_2$ be a character whose conductor is $d_2$ with $d_2 |d$. Then
$$\sum_{y \in (\mathbb Z/d \mathbb Z)^\times} \overline \theta_2 (y) e_{\theta_1}
\log_p(\alpha^y - \zeta_{p^{n+1}})=
\left( \sum_{y \in (\mathbb Z/d_2 \mathbb Z)^\times} \overline \theta_2 (y) e_{\theta_1}
\log_p(\alpha_2^y - \zeta_{p^{n+1}}) \right) x(\theta)$$ where $\alpha_2=\alpha^{d/d_2}$ and $x(\theta) \in \zp[\theta]$.
\end{lemma}
\begin{proof}It is sufficient to consider the case where $d=ld_2$ with $l$ a prime number.
Let $S$ be the left-hand side sum in lemma 1.8.
\begin{itemize}
\item First case : $l\not|d_2$.  For $y \in (\mathbb Z/d_2 \mathbb Z)^\times$
we want to write the elements $z$ of $(\mathbb Z/d\mathbb
Z)^\times$ as $z=y+k d_2$ for $0\leq k\leq l-1$. Howerver for all
$y$ there exists a $k_y, \, 0 \leq k_y \leq l-1$, such that $y+k_y
d_2 \equiv 0 \mod l$ (i.e. $y+k_y d_2 \not\in (\mathbb Z/d\mathbb
Z)^\times$). We can write $S$ as $$S=\sum_{y \in (\mathbb
Z/d_2\mathbb Z)^\times}\overline{\theta}_2(y) e_{\theta_1} \log_p
\prod_{\stackrel{k=0}{k\neq
k_y}}^{l-1}(\alpha^{kd_2+y}-\zeta_{p^{n+1}}).$$ Note that
$\alpha^{d_2}$ is a $l$th root of unity. Since $\prod_{\zeta \in
\mu_l}(X-\zeta Y)=(X^l-Y^l)$ we get $S$
\begin{eqnarray*}S&=&\sum_{y \in (\mathbb Z/d_2\mathbb Z)^\times}\overline{\theta}_2(y) e_{\theta_1} \log_p
\left( \frac{\alpha_2^{y}-\zeta_{p^{n+1}}^l}{\alpha_2^{\frac{y-k_y d_2}l}-\zeta_{p^{n+1}}}\right)\\
&=&(\theta_1(l)-\overline{\theta_2}(\sigma_{\frac{y+k_y d_2}l}))
\sum_{y \in (\mathbb Z/d_2\mathbb Z)^\times}\overline{\theta}_2(y) e_{\theta_1} \log_p (\alpha_2^y-\zeta_{p^{n+1}}).
\end{eqnarray*}
\item Second case : $l|d_2$. Then all $z=y+k d_2$ with $y \in (\mathbb Z/d_2\mathbb Z)^\times$ and $0 \leq k \leq l-1$
are invertible modulo $d$. Thus
\begin{eqnarray*}
S&=& \sum_{y \in (\mathbb Z/d_2\mathbb Z)^\times}\overline{\theta}_2(y) e_{\theta_1} \log_p(\alpha_2^y-\zeta_{p^{n+1}}^l)\\
&=& \theta_1(l) \sum_{y \in (\mathbb Z/d_2\mathbb Z)^\times}\overline{\theta}_2(y) e_{\theta_1} \log_p (\alpha_2^y-\zeta_{p^{n+1}}).
\end{eqnarray*}
\end{itemize}
\end{proof}

For a given $d$ dividing $p-1$ the two previous lemmas yield the
equality
\begin{equation}\tau(\overline \theta_2) x(\theta_1 \theta_2)\epsilon_n(\theta_1 \theta_2) e_{\theta_1} \mathscr T_n =
-\sum_{y\in (\mathbb Z/d \mathbb Z)^\times} \overline{\theta }_2(y) e_{\theta_1}
\log_p(\alpha^y - \zeta_{p^{n+1}}).
\end{equation}

\begin{lemma} Let $d$ be an integer dividing  $p-1$ and $\alpha$ a primitive $d$th root of unity.
Fix $\theta_1 \in \widehat \Delta$ different from $1, \omega$.
Then there exists $u_n \in \zp[\Gamma_n]$ such that $$u_n
e_{\theta_1} \mathscr T_n= e_{\theta_1} \log_p(\alpha-
\zeta_{p^{n+1}}).$$
\end{lemma}
\begin{proof}
Suppose first that $\theta_1$ is even. Summing the equality (2)
over all the $\theta_2$ such that the product $\theta_1 \theta_2$
is even we obtain $$\left(
\sum_{\stackrel{f_{\theta_2}|d}{\theta_2 \textrm{ even
}}}\tau(\overline \theta_2) x(\theta_1
\theta_2)\epsilon_n(\theta_1 \theta_2) \right) e_{\theta_1}
\mathscr T_n = -\sum_{y\in (\mathbb Z/d \mathbb Z)^\times} \left(
\sum_{\stackrel{f_{\theta_2}|d}{\theta_2 \textrm{ even
}}}\overline{\theta }_2(y) \right) e_{\theta_1} \log_p(\alpha^y -
\zeta_{p^{n+1}}).$$ We have $$
\sum_{\stackrel{f_{\theta_2}|d}{\theta_2 \textrm{ even
}}}\overline{\theta }_2(y) = \left\{
\begin{array}{ll} 0 & \textrm{when $y \not\equiv \pm 1 \mod d$}\\
\frac{\varphi(d)}{2} & \textrm{when $y \equiv \pm 1 \mod d$}.
\end{array} \right.$$ Hence the right-hand side of the above is equal to
$-\frac{\varphi(d)}{2} e_{\theta_1} \log_p \left[(\alpha-
\zeta_{p^{n+1}})(\alpha^{-1}-\zeta_{p^{n+1}}) \right]$. Since
$\theta_1$ is even we have
$e_{\theta_1}\log_p(\alpha^{-1}-\zeta_{p^{n+1}})=e_{\theta_1}\log_p(\alpha-\zeta_{p^{n+1}}^{-1})=e_{\theta_1}\log_p(\alpha-\zeta_{p^{n+1}})$.
The right-hand side thus equals $-\varphi(d)e_{\theta_1}
\log_p(\alpha-\zeta_{p^{n+1}})$. Set $$u_n=-\frac 1{\varphi(d)}
\left(\sum_{\stackrel{f_{\theta_2}|d}{\theta_2 \textrm{ even
}}}\tau(\overline \theta_2) x(\theta_1
\theta_2)\epsilon_n(\theta_1 \theta_2) \right).$$ As $u_n$ is
Galois-invariant we have $u_n \in
\zp[\Gamma_n]$. The result follows in this case. \\ \indent When
$\theta_1$ is odd, the proof runs the same except that $$
\sum_{\stackrel{f_{\theta_2}|d}{\theta_2 \textrm{ odd
}}}\overline{\theta }_2(y) = \left\{ \begin{array}{ll} 0 &
\textrm{if $y \not\equiv \pm 1 \mod d$}\\ \pm\frac{\varphi(d)}{2}
& \textrm{if $y \equiv \pm 1 \mod d$} \end{array} \right.$$ and
$e_{\theta_1}\log_p(\alpha^{-1}-\zeta_{p^{n+1}})=-e_{\theta_1}\log_p(\alpha-\zeta_{p^{n+1}}^{-1})$.
\end{proof}

\begin{theorem} Let $\theta \in \widehat \Delta$, $\theta \neq 1, \omega$.
Let $U_n$ be the group of principal units of $\mathbb Q
_p(\zeta_{p^{n+1}})$ defined by $U_n=1+(\zeta_{p^{n+1}}-1)
\zp[\zeta_{p^{n+1}}]$ and $\mathscr T_n = \sum_{i=0}^n p^{i-n}
\zeta_{p^{i+1}}$. Then for all $n \geq 0$, $$e_\theta \log_p U_n =
\zp[\Gamma_n] e_\theta \mathscr T_n.$$
\end{theorem}
\begin{proof}By lemma 1.9 there exists $u_n \in \zp[\Gamma_n]$ such that $u_n e_\theta \mathscr T_n =
 e_\theta \log_p(\alpha-\zeta_{p^{n+1}})$. By the results of \cite{WAS} section 13.8
there exists an integer $d$ dividing $p-1$ and a primitive $d$th
root of unity $\alpha$ such that $e_\theta \log_p U_n =
\zp[\Gamma_n] e_\theta \log_p(\alpha-\zeta_{p^{n+1}})$. Hence
$e_\theta \log_p U_n \subseteq \zp[\Gamma_n] e_\theta \mathscr
T_n$. The converse inclusion will be proved in the following
section (lemma 1.11 and proposition 1.12). Just note it implies
$u_n \in \zp[\Gamma_n]^\times$ for such an $\alpha$.

\end{proof}

We now derive some corollaries from theorem 1.10, the first of
which is a well-known result of Iwasawa.
\begin{Cor}[Iwasawa](\cite{IWA2} or \cite{WAS} theorem 13.56) Let $\theta$ be an even nontrivial character in $\widehat \Delta$. Let
$(\epsilon_n(\theta))_{n\geq 0} \in \limproj \zp[\Gamma_n]$ be the
element corresponding to the power series $f(\frac{1+p}{1+T}~-~1~,
\theta)$. Let $\overline C_n$ the closure of the cyclotomic units
in $U_n$. For all $n \geq 0$ we have an isomorphism $$e_\theta
U_n/\overline C_n \simeq \zp[\Gamma_n]/\epsilon_n(\theta)
\zp[\Gamma_n].$$
\end{Cor}
\begin{proof}By basic results on cyclotomic units we have
$$e_\theta\log_p \overline C_n = \zp[\Gamma_n]e_\theta
\log_p(1-\zeta_{p^{n+1}}).$$ Lemma 1.9 allpied to $d=1$ shows that
there exists $u_n \in \zp[\Gamma_n]$ such that $u_n e_\theta
\mathscr T_n = e_\theta \log_p(1-\zeta_{p^{n+1}})$. The definition
of $u_n$ gives $u_n=-\epsilon_n(\theta)$. We deduce that $\log_p
\overline C_n = \zp[\Gamma_n]\epsilon_n(\theta) e_\theta \mathscr
T_n$. Moreover theorem 2 shows that $\log_p U_n = \zp[\Gamma_n]
e_\theta \mathscr T_n$. As $\theta$ is even and nontrivial we have
isomorphisms $$e_\theta U_n/\overline C_n \simeq e_\theta \log_p
U_n/\log_p \overline C_n \simeq \zp[\Gamma_n]/\epsilon_n(\theta)
\zp[\Gamma_n].$$
\end{proof}

\begin{Cor} Assume $d \geq 3$. There exists an odd character $\theta_2$ whose conductor divides $d$ such that
there are at least $\sqrt p -2$ odd character $\theta_1 \in
\widehat \Delta$, $\theta_1 \neq \omega$ satisfying $$ f(T,
\theta_1\theta_2) \in W \entL T \entR ^\times.$$
\end{Cor}

\begin{proof} By a result of Angl\`es (\cite{Ang}, theorem 5.4) there are at least $\sqrt p -2$
odd characters $\theta_1 \in \widehat \Delta$, such that $\theta
\neq \omega$ and for all $n\geq 0$, $$e_{\theta_1} \log_p U_n =
e_{\theta_1} \log_p(\alpha-\zeta_{p^{n+1}}).$$ For such a
character we know that $$u_n =
\sum_{\stackrel{f_{\theta_2}|d}{\theta_2 \textrm{ odd
}}}\tau(\overline \theta_2) x(\theta_1
\theta_2)\epsilon_n(\theta_1 \theta_2) \in \zp[\Gamma_n]^\times
.$$ Therefore there exists at least one $\theta_2$ such that
$\epsilon_0(\theta_1 \theta_2) \in W^\times$. Since $\theta_1
\theta_2$ is a character of the first kind it follows that $f(T,
\theta_1\theta_2) \in W \entL T \entR ^\times$.
\end{proof}
A similar proof yields the following result.
\begin{Cor} Let $\theta_1 \in \widehat \Delta$, $\theta_1 \neq 1, \omega$.
There exists $\theta_2$ whose conductor divides $p-1$ such that
$\theta_1 \theta_2$ is even and the generalized Bernoulli number
$B_{1, \theta_1 \theta_2 \omega^{-1}}$ is rpime to $p$
\end{Cor}

\subsection{Some index computations}
Let us recall some notations. As usual $p$ is an odd prime number
and $\on$ is the ring of integers of $\mathbb Q_p(\zeta_p^{n+1})$.
Let $\mathscr T_n$ be $\sum_{i=0}^n p^{i-n} \zeta_{p^{i+1}}$,
$e_0=\frac{1}{p^n}\sum_{\gamma \in \Gamma_n} \gamma$ and for $1
\leq d \leq n$ $$e_d=\sum_{\stackrel{\chi \in
\widehat{\Gamma}_n}{\chi \textrm{ of conductor }p^{d+1}}}
e_\chi.$$ We first want to compute the index
$[\frac{1}{p^n}e_\theta \on : \zp[\Gamma_n] e_\theta \mathscr
T_n]$ with $\theta \in \widehat{\Delta}$. Let $\Lambda_n$ the
maximal order of $\qp[\Gamma_n]$. Then $\Lambda_n[\Delta]$ is the
maximal order of $\qp[\Delta \times \Gamma_n]$ and the Leopoldt
theorem (see \cite{Erez}) states that $\on$ is a free
$\Lambda_n[\Delta]$-module generated by $T_n=\sum_{i=0}^n
\zeta_{p^{i+1}}$.

We can therefore substitute $\frac{1}{p^n}e_\theta \on$ by
$\Lambda_n e_\theta \sum_{i=0}^n \frac{\zeta_{p^{i+1}}}{p^n}$ in
our calculations. The index $[\frac{1}{p^n}e_\theta \on :
\zp[\Gamma_n] e_\theta \mathscr T_n]$ with $\theta \in
\widehat{\Delta}$ is the product $$[\Lambda_n e_\theta
\sum_{i=0}^n \frac{\zeta_{p^{i+1}}}{p^n} : \Lambda_n e_\theta
\sum_{i=0}^n \frac{\zeta_{p^{i+1}}}{p^{n-i}}] . [\Lambda_n
e_\theta \mathscr T_n :\zp[\Gamma_n] e_\theta \mathscr T_n].$$ By
a standard calculation of discriminants we have $[\Lambda_n
e_\theta \mathscr T_n :\zp[\Gamma_n] e_\theta \mathscr
T_n]=[\Lambda_n : \zp[\Gamma_n]]=p^{\frac{p^n-1}{p-1}}$. In order
to find $[\Lambda_n e_\theta \sum_{i=0}^n
\frac{\zeta_{p^{i+1}}}{p^n} : \Lambda_n e_\theta \sum_{i=0}^n
\frac{\zeta_{p^{i+1}}}{p^{n-i}}]$ we notice that
 $\Lambda_n=\bigoplus_{d=0}^n \zp[\gamma_n]e_d$. Thus
\begin{eqnarray*}
[\Lambda_n e_\theta \sum_{i=0}^n \frac{\zeta_{p^{i+1}}}{p^n} :
\Lambda_n e_\theta \sum_{i=0}^n \frac{\zeta_{p^{i+1}}}{p^{n-i}}] &=&
\prod_{i=0}^n [\zp[\Gamma_n] e_i
\frac{\zeta_{p^{i+1}}}{p^n} :\zp[\Gamma_n] e_i
\frac{\zeta_{p^{i+1}}}{p^{n-i}}] \\
 &=& \prod_i p^{i(p^i-p^{i-1})} ,
\end{eqnarray*}
because the $\zp$-rank of $\zp[\Gamma_n] \zeta_{p^{i+1}}$ is
$p^i-p^{i-1}$. Moreover we have the equality \\$\sum_{i=0}^n
d(p^i-p^{i-1})=n p^n-\frac{p^n-1}{p-1}$.

We have thus shown the following result.
\begin{lemma} Let $\theta \in \widehat{\Delta}$. We have
$[\frac{1}{p^n}e_\theta \on : \zp[\Gamma_n] e_\theta
\mathscr T_n]=p^{np^n}$.
\end{lemma}

The next index we want to calculate is $[\frac{1}{p^n} e_\theta
\on : e_\theta \log_p U_n]$.
\begin{prop} Let $\theta \in \widehat{\Delta}$. We have
\begin{displaymath}[\frac{1}{p^n} e_\theta \on : e_\theta \log_p(U_n)] = \left\{
\begin{array}{ll} p^{np^n} & \textrm{when $\theta \neq 1, \omega$} \\
p^{np^n+1} & \textrm{when $\theta=1$}\\ p^{np^n+n+1} &\textrm{
when $\theta=\omega$}.
\end{array} \right. \end{displaymath}
\end{prop}

\begin{proof}
The proof is based on a result of John Coates, see \cite{FRO}. We
first show the inclusion $\log_p U_n \subseteq \frac{1}{p^n} \on$.
Define $\pi_n= \zeta_{p^{n+1}}-1$ and let $u \in U_n$. Then there
exists $a \in \mathbb Z$ such that  $$u \zeta_{p^{n+1}} = 1 \mod
\pi_n^2.$$ Therefore $\log_p U_n = \log_p(1+\pi_n^2 \on)$.
Moreover  we have $u \equiv 1 \mod\pi_n^2$ implies $v_p(u^{p^n}-1)
\geq \frac{2}{p-1}$ for $u \in U_n$. Then by lemma 5.5 of
\cite{WAS} we have $\log_p U_n \subseteq \frac{1}{p^n} \on$.

We now calculte the index. We introduce the group $V=1+p \on$. By
standard properties of the $p$-adic logarithm we know that $\log_p
V=p \on$. Then $$[\frac{1}{p^n} e_\theta \on : e_\theta \log_p
V]=[\frac{1}{p^n} e_\theta \on : p e_\theta \on]= p^{p^n (n+1)}.$$
It remains to compute $[e_\theta \log_p U_n : e_\theta \log_p V]$.
Consider the morphism $$ (1+ \pi_n^2 \on)/V
\stackrel{\log_p}{\longrightarrow} \log_p U_n /\log_p V,$$ the
kernel of which consists of the roots of unity. Then for $\theta
\neq \omega$ we have $$e_\theta \log_p U_n /\log_p V \simeq
e_\theta (1+ \pi_n^2 \on)/V = \prod_{i=2}^{p^n(p-1)-1} e_\theta
\frac{1+\pi_n^i \on}{1+\pi_n^{i+1} \on}.$$

For all integers $i \geq 1$ we have an isomorphism of
$\zp$-modules $$\begin{array}{ccl} \frac{1+\pi_n^i
\on}{1+\pi_n^{i+1} \on} & \longrightarrow & \on/ \pi_n \on \simeq
\mathbb F_p \\ 1+x \pi_n^i & \mapsto &  x \mod \pi_n
.\end{array}$$ Hence the $\theta$-part of $\frac{1+\pi_n^i
\on}{1+\pi_n^{i+1} \on}$ is $\mathbb F_p$ or $\{0\}$. Moreover for
all $\delta \in \Delta$ we have $x^\delta \equiv x \mod \pi_n$.
Thus for $\theta \in \widehat{Delta}$ we have $$
(1+x\pi_n^i)^{e_\theta} \equiv \prod_{\delta \in \Delta}
(1+\frac{1}{p-1} \overline{\theta}(\delta)x(\pi_n^i)^\delta)
\mod(\pi_n^{i+1}). $$ However $(\pi_n^i)^{\delta} \equiv
\omega^i(\delta)\pi_n^i \mod(\pi_n^{i+1})$ which yields
$$(1+x\pi_n^i)^{e_\theta} \equiv 1-(\sum_{\delta\in \Delta}
\overline{\theta}(\delta)\omega^i(\delta)) x \pi_n^i
\mod(\pi_n^{i+1}).$$ Thus for $\theta=\omega^k$ we get
\begin{displaymath}e_\theta \frac{1+\pi_n^i \on}{1+\pi_n^{i+1} \on} = \left\{
\begin{array}{ll} 0 & \textrm{when $i \not\equiv k \mod (p-1)$} \\
\mathbb F_p & \textrm{when $i \equiv k \mod (p-1)$},\\
\end{array} \right. \end{displaymath}

and $[e_\theta (1+\pi_n^2 \on) : e_\theta (1+p \on)] = p^{C(k)}$
where $C(k)$ is the number of integers $i$ such that $2\leq i \leq
p^n (p-1)-1$ and $i \equiv k \mod(p-1)$. When $k\neq 0,1$ (i.e.
$\theta \neq \omega, 1$) we have $C(k)=p^n$ and when $k=0,1$ we
have $C(k)=p^n-1$. Hence the result for $\theta \neq \omega$.

The case $\theta = \omega$  is similar except that $\zeta_p,
\cdots,\zeta_{p^n}$ are in the kernel of the morphism $\log_p :
(1+ \pi_n^2 \on)/V \longrightarrow \log_p U_n /\log_p V $.
\end{proof}
This completes the proof of theorem 1.10. We might want to
reformulate it according to Leopoldt's element
$T_n=\sum_{i=0}^n\zeta_{p^{i+1}}$. Let $l_n=\sum_{i=0}^n
p^{n-i}e_d \in \mathbb Q_p[\Gamma_n]$. This element is constructed
to satisfy the identity $l_n \mathscr T_n=T_n$. Notice that $l_n
\circ (\sum_{i=0}^n p^{i-n}e_d)=1$ so that $l_n \in \mathbb
Q_p[\Gamma_n]^\times$. Let $\mathscr L_n =l_n \circ \log_p$.
\begin{Cor} Let $\theta \in \widehat{\Delta}$with $\theta \neq 1, \omega$.
Then, $e_\theta \mathscr L_n U_n = \zp[\Gamma_n] e_\theta T_n$.
\end{Cor}

\subsection{The case of the Teichmüller character}
For technical reasons the results in this section are only valid
when $p \geq 5$. Let us begin with the following proposition.
\begin{prop} There exists $\alpha \in \mu_{p-1} \setminus \{\pm1\}$ such that
for all $n \geq 0$,
$$e_\omega \log_p U_n = \zp[\Gamma_n] e_\omega \log_p (\alpha -\zeta_{p^{n+1}}).$$
\end{prop}
\begin{proof} A careful reading of the proof of theorem 13.54 in
\cite{WAS} shows that $\limproj_{n \geq 0} e_\omega \log_p U_n$ is
a free $\Lambda$-module of rank $1$ and that when
$(\epsilon_n)_{n\geq 0} \in  \limproj_{n \geq 0} e_\omega \log_p
U_n$ is such that $e_\omega \log_p U_0 = \zp[\epsilon_0]$ then for
all $n\geq 0$, $$e_\omega \log_p U_n = \zp[\Gamma_n] \epsilon_n.$$
Therefore it is sufficient to show there exists $\alpha \in
\mu_{p-1}$, $\alpha \neq \pm1$ such that
\begin{equation}e_\omega \log_p U_0 = \zp e_\omega \log_p(\alpha - \zeta_p).
\end{equation}
Let us prove $(3)$. Recall that $e_\omega \log_p U_0 = e_\omega
\pi_n^2 \mathcal O_0=p\,\tau(\omega^{-1}) \zp$. Let $\pi$ denote
the element $\zeta_p-1$.
\begin{lemma} Let $x \in \pi^2 \zp[\zeta_p]$. Then $\log_p(1+x) \equiv \frac{(1+x)^p-1}{p} \mod (p\pi^2)$.
\end{lemma}
\begin{proof}
We have the congruences
\begin{itemize}
\item for all $n \geq p$, $v_\pi(\frac{x^n}{n})\geq p+1$ and
\item for $1\leq k\leq p-1$, $\frac{\binom{p}{k}}{p} \equiv \frac{(-1)^{k+1}}{k} \mod p.$
\end{itemize}
We then have $$\begin{array}{cl} \log_p(1+x) & = \sum_{n \geq1}
\frac{(-1)^{n+1}}{n}x^n \\
 & \equiv \sum_{n=1}^{p-1} \frac{\binom{p}{n}}{p} x^n \mod(p\pi^2) \\
 & \equiv \frac{(1+x)^p-1-x^p}{p} \mod(p\pi^2)\\
 & \equiv \frac{(1+x)^p-1}{p} \mod(p \pi^2) ,
\end{array}$$
and the lemma is proved.
\end{proof}
\begin{lemma}Let $\alpha \in \mu_{p-1} \setminus \{1\}$ and
$\theta \in \widehat{\Delta}$, $\theta \neq 1$. We have $$
e_\theta \log_p(\alpha-\zeta_p) \equiv
-\frac{\theta(-1)}{\omega(\alpha-1)} \tau(\overline{\theta})\left(
\sum_{k=1}^{p-1} (-1)^k \frac{\binom{p}{k}}{p} \theta(k) \alpha^k
\right) \mod(p\pi^2).$$
\end{lemma}
\begin{proof}
Let $\gamma = \frac{\alpha-\zeta_p}{\omega(\alpha-1)}
\zeta_p^{\overline{\omega}(\alpha-1)}$. Notice that
$\log_p(\gamma)=\log_p(\alpha - \zeta_p)$ and $\gamma \equiv 1
\mod(\pi^2)$, which allows us to apply lemma 1.14. We get
$$\log_p(\alpha - \zeta_p) \equiv
\frac{\frac{(\alpha-\zeta_p)^p}{\omega(\alpha-1)} -1}{p} \mod
(p\pi^2).$$ Taking the $\theta$-part and expanding the sum yields
the required result.
\end{proof}
Since $e_\omega  \log_p U_0 = p \,\tau(\omega^{-1}) \zp$ and
$v_\pi(\tau(\omega^{-1}))=1$ we deduce from lemma 1.15 the
following result.
\begin{lemma}Let $\alpha \in \mu_{p-1} \setminus \{1\}$. The element $e_\omega \log_p(\alpha - \zeta_p)$
is a generator of $e_\omega \log_p U_0$ if and only if
$$\sum_{k=1}^{p-1}(-1)^k \frac{\binom{p}{k}}{p} \omega(k) \alpha^k \not\equiv 0 \mod (p^2).$$
\end{lemma}

\begin{lemma}There exists $\alpha \in \mu_{p-1} \setminus \{\pm1\}$ such that
$\sum_{k=1}^{p-1}(-1)^k \frac{\binom{p}{k}}{p} \omega(k) \alpha^k \not\equiv 0 \mod (p^2)$.
\end{lemma}

This is precisly what we need to complete the proof of proposition
1.13.
\begin{proof}
Let us consider the polynomial $P(X)=\sum_{k=1}^{p-1} (-1)^k
\frac{\binom{p}{k}}p \omega(k) X^k \in \zp[X]$. We know that
$\binom{p}{k}/p \equiv (-1)^{k+1}/k \mod p$ and thus $$P(X) \equiv
-X \prod_{\alpha \in \mu_{p-1} \setminus \{1\}} (X-\alpha) \mod
p.$$ We can therefore use Hensel's lemma which shows that $P(X)=-X
\prod_{\alpha \in \mu_{p-1} \setminus \{1\}}(X-a(\alpha))$ where
$a(\alpha) \in \zp^\times$ and $a(\alpha) \equiv \alpha \mod p$.
Notice that $P(-1) \equiv 0 \mod(p^2)$.

Let us assume the lemma is false. Then for all $\alpha \in
\mu_{p-1} \setminus \{1\}$, we have $a(\alpha) \equiv \alpha \mod
p^2$. Thus $P(X) \equiv -X \frac{X^{p-1}-1}{X-1} \mod p^2$.
Comparing with the expression of $P(X)$ we deduce that for all $k
\in \{1, \cdots,p-1\}$, we have $$(-1)^k \frac{\binom{p}{k}}{p}
\omega(k) \equiv -1 \mod p^2.$$

Let us apply this congruence with $k=2$ and $k=4$. The assumption
$p \geq 5$ is essenial in what follows. We obtain on the one hand
$(p-1) \omega(2)/2 \equiv -1 \mod p^2$, that is to say
$\omega(2)/2 \equiv p+1 \mod p^2$. On the other hand we have
$\omega(4)/4 \equiv 1+(11/6)p \mod p^2$. But
$\omega(4)/4=(\omega(2)/2)^2\equiv 1+2p \mod p^2$. This implies
$11/6 \equiv 2 \mod p^2$ which is not possible.
\end{proof}
This finishes the proof of proposition 1.13.
\end{proof}
In order to apply some of the results of section 1.2, we let
$\theta_2$ be an odd character of conductor $d$ with $d$ dividing
$p-1$. Recall that $$\tau(\overline{\theta_2}) \epsilon_n(\theta_2
\omega) e_\omega \mathscr T_n = - \sum_{y \in (\mathbb Z/ d
\mathbb Z)^\times} \overline{\theta_2}(y) e_\omega \log_p(\alpha^y
- \zeta_{p^{n+1}})$$ where $\alpha$ is a primitive $d$th root of
unity and $\mathscr T_n = \sum_{i=0}^n p^{i-n} \zeta_{p^{i+1}}$.

Since $\theta_2(p)=1$ we have for all characters $\chi$ and all
integers $n \geq 1$, $$L_p(1-n, \chi)=-(1-\chi
\omega^{-n}(p)p^{n-1})\frac{B_{n,\chi \omega^{-n}}}{n}$$ from
which we deduce immediatly that $L_p(0, \theta_2 \omega)=0$. Thus
if $f(T, \theta_2 \omega)$ denotes the power series associated to
$L_p(s, \theta_2\omega)$ we have the factorization $f(T, \theta_2
\omega)=T h(T, \theta_2 \omega)$. Then there exists $H(T, \theta_2
\omega) \in W\entL T \entR$ such that $$f(\frac{1+pd}{1+T}-1,
\theta_2 \omega)=(T-pd) H(T, \theta_2 \omega).$$ As usual $T$
corresponds to $\sigma_{1+pd}-1 \in W\entL
\Gal(K_\infty/K_0)\entR$. Let $\widetilde{\epsilon_n}$ be the
element of $\limproj W[\Gamma_n]$ associated to $H$ via the
isomorphism in lemma 1.1. We have
\begin{eqnarray}
\epsilon_n(\theta_2 \omega) e_\omega \mathscr T_n &=&
\widetilde{\epsilon_n}(\theta_2 \omega) (\sigma_{1+pd}-1-pd)
e_\omega \mathscr T_n \\ &=& \widetilde{\epsilon_n}(\theta_2
\omega)(\sigma_{1+p}^s-(1+p)^s) e_\omega \mathscr T_n ,
\end{eqnarray}
where $s=\log_p(1+pd)/\log_p(1+p) \in \zp^\times$. So we have
$$\sigma_{1+p}^s-(1+p)^s=(\sigma_{1+p}-1-p).u$$ with $u\in
\zp[\Gamma_n]^\times$. The same kind of calculation than in the general case
shows that there exists a unit $ u_n \in \zp[\Gamma_n]$ such that
$$u_n (\gamma_0-1-p) e_\omega \mathscr{T}_n = e_\omega
\log_p(\alpha-\zeta_{p^{n+1}})$$ where $\gamma_0 = \sigma_{1+p}$
and $\alpha$ is a primitive $d$th root of unity.

\begin{lemma} We have
$$[\zp[\Gamma_n] e_\omega \mathscr T_n : \zp[\Gamma_n] (\gamma_0-1-p) e_\omega \mathscr T_n]=
p^{np^n+n+1}.$$
\end{lemma}
\begin{proof}
By lemma 1.11 we have $[\frac{1}{p^n}e_\omega \on : \zp[\Gamma_n]
e_\omega \mathscr T_n]=p^{np^n}$. Also
\begin{eqnarray*}
[\zp[\Gamma_n] e_\omega \mathscr T_n : \zp[\Gamma_n] (\gamma_0-1-p) e_\omega \mathscr T_n] &=&
[\Lambda/\omega_n \Lambda : (T-p) \Lambda/\omega_n \Lambda] \\ &=& [\Lambda : (\omega_n, T-p)]
\end{eqnarray*}
where $\Lambda \simeq \zp \entL T\entR$ is the Iwasawa algebra.
The $\Lambda$-module $M=\Lambda/(T-p)$ has no $\zp$-torsion. Its
characteristic polynomial is $T-p$ and it is a standard result
that $|M/\omega_n M|=\prod_{\zeta \in \mu_{p^n}}(\zeta-p-1)$. Then
$$[\Lambda : (\omega_n, T-p)] \sim p \prod_{\stackrel{\zeta \in
\mu_{p^n}}{\zeta\neq 1}} (\zeta-1) \sim p^{n+1}$$ where $\sim$
means 'has same $p$-adic valuation as'.

\end{proof}
We have thus established the following result.
\begin{theorem} Let $p \geq 5$ be a prime number. We have
$$e_\omega \log_p U_n = \zp[\Gamma_n] (\gamma_0-1-p) e_\omega
\mathscr T_n$$ where $\mathscr T_n = \sum_{i=0}^n p^{i-n}
\zeta_{p^{i+1}}$. With the notations ofsection 1.4, this writes as
$$e_\omega \mathscr L_n U_n = \zp[\Gamma_n] (\gamma_0-1-p)
e_\omega T_n$$ where $T_n=\sum_{i=0}^n \zeta_{p^{i+1}}$.
\end{theorem}

\subsection{The case of the trivial character}
The main difference between the trivial character and the other
ones is that the power series associated to $L_p$ has not integral
coefficients. Thus We have to work with the power series $g(T)$
defined by $g(T)=(1-\frac{1+q_0}{1+T})f(T,1) \in \Lambda$. See
\cite{WAS}, chapter 7 for more details on this particularly
proposition 7.9. In our situation we have $q_0=p$. Let
$h(T)=g(\frac{1+p}{1+T}-1)=-T f(\frac{1+p}{1+T}-1,1)$ and $\eta_n
\in \limproj \zp[\Gamma_n]$ be the element corresponding to $h$.
This means that we have, as in section 1.1, for all $\chi \in
\widehat{\Gamma}_n$ different from $1$, $$e_\chi \eta_n =
(1-\chi(\gamma_0)) L_p(1,\chi) e_\chi.$$ Let
$T_\Delta=\sum_{\delta \in \Delta} \delta$ the norm element of
$\zp[\Gal(K_0)/\mathbb Q]$.
\begin{prop}Let $E_n$ be the group of units in $K_n$, $\overline{E}_n$ its closure in $U_n$ and
$ \widetilde{\mathscr T}_n=\sum_{i=1}^n p^{i-n} \zeta_{p^{i+1}}$ (note that the sum begins at $i=1$).
We have $$T_\Delta \log_p \overline{E}_n=\zp[\Gamma_n]T_\Delta \widetilde{\mathscr T}_n.$$
\end{prop}
\begin{proof}
Let the field $\mathbb B_n \in \mathbb Q(\zeta_{p^{n+1}})$ be such
that $[\mathbb B_n : \mathbb Q]=p^n$. Iwasawa has showed
\cite{IWA} that $p$ does not divide the class number of $\mathbb
B_n$. Thus \cite{WAS}, theorem 8.2 shows that $$ T_\Delta \log_p
\overline{E}_n = T_\Delta \log_p \overline{C}_n$$ where
$\overline{C}_n$ is the closure of $C_n$ in $U_n$.

Moreover we have $T_\Delta \log_p \overline C_n = \zp[\Gamma_n]
T_\Delta \log_p((1-\zeta_{p^{n+1}})^{\gamma_0-1})$. Now fix a
character $\chi \in \widehat{\Gamma}_n$ different from $1$ whose
conductor is $f_\chi=p^{k+1}$ with $1 \leq k\leq n$. We have
$$\frac{1}{p^n} \tau(\overline{\chi})=e_\chi T_\Delta
\widetilde{\mathscr T}_n$$ and then
$$\frac{1}{p^n}(1-\chi(\gamma_0))\tau(\overline
\chi)L_p(1,\chi)=e_\chi \eta_n T_\Delta \widetilde{\mathscr
T}_n.$$ A little calculation gives
\begin{eqnarray*}\tau(\overline \chi)L_p(1,\chi)&=&-\sum_{\stackrel{\delta \in \Delta}{\gamma\in \Gamma_d}}
\overline \chi(\delta) \log_p(1-\zeta_{p^{d+1}}^{\delta \gamma})\\
&=&-p^n e_\chi T_\Delta \log_p(1-\zeta_{p^{n+1}}).
\end{eqnarray*}
This proves that for all $\chi \in \widehat{\Gamma}_n$, $\chi \neq
1$, we have the following equality $$e_\chi T_\Delta
\log_p((1-\zeta_{p^{n+1}})^{\gamma_0-1})= -e_\chi T_\Delta
\widetilde{\mathscr T}_n.$$ We check that this equality also holds
when $\chi =1$ and summing over all the characters gives the
required result.
\end{proof}

\begin{lemma}We have $$[T_\Delta\frac{1}{p^n}\on : \zp[\Gamma_n] T_\Delta(p \zeta_p + \widetilde{\mathscr T}_n)]=
p^{np^n+n+1}.$$
\end{lemma}
\begin{proof}We already know that $[\Lambda_n : \zp[\Gamma_n]]=p^{\frac{p^n-1}{p-1}}$
where $\Lambda_n$ is the maximal order of $\qp[\Gamma_n]$. It
remains to calculate $[\frac 1 {p^n}T_\Delta \on : \Lambda_n
T_\Delta(p \zeta_p + \widetilde{\mathscr T}_n)]$. Decompose into
characters yields $$[\frac 1 {p^n}T_\Delta \on : \Lambda_n
T_\Delta(p \zeta_p + \widetilde{\mathscr T}_n)] = \prod_{d=0}^n
[e_d\frac 1 {p^n}T_\Delta \on : e_d\Lambda_n T_\Delta(p \zeta_p +
\widetilde{\mathscr T}_n)].$$ For $1 \leq d \leq n$ we have $e_d
\frac 1 {p^n}T_\Delta \on = \zp[\Gamma_n]T_\Delta(p^{-n}
\zeta_{p^{d+1}})$ and $e_d\Lambda_n T_\Delta(p \zeta_p +
\widetilde{\mathscr T}_n)=\zp[\Gamma_n]T_\Delta(p^{d-n}
\zeta_{p^{d+1}})$. Then when $d \neq 0$ we have $[e_d\frac 1
{p^n}T_\Delta \on : e_d\Lambda_n T_\Delta(p \zeta_p +
\widetilde{\mathscr T}_n)]=p^{d(p^d-p^{d-1})}$. For $e_0$ we
notice that $e_0T_\Delta \frac 1 {p^n} \on=\frac 1{p^n}\zp$ and
$e_0\widetilde{\mathscr T}_n=0$. Then we have $[e_0\frac 1
{p^n}T_\Delta \on : e_0\Lambda_n T_\Delta(p \zeta_p +
\widetilde{\mathscr T}_n)]=[\frac 1{p^n}\on : p \zp]=p^{n+1}$. The
lemma is now proved.
\end{proof}
Recall that proposition 1.12 gives the index $$[\frac 1{p^n}
T_\Delta \on : T_\Delta \log_p U_n]=p^{np^n+1}.$$ Moreover note
that $p\zeta_p \in \log_p U_n$ so that $T_\Delta(p \zeta_p +
\widetilde{\mathscr T}_n) \in T_\Delta \log_p U_n$. However the
value of the index furnished by proposition 1.12 implies that the
$\zp[\Gamma_n]$-module generated by $T_\Delta(p \zeta_p +
\widetilde{\mathscr T}_n) $ cannot be equal to $T_\Delta \log_p
U_n$. Let us define $\mathscr M = \zp[\Gamma_n] T_\Delta p\zeta_p
+ \zp[\Gamma_n]T_\Delta \widetilde{\mathscr T}_n$. We check that
$\mathscr M \subseteq T_\Delta \log_p U_n$.

Now we want to calculate the index $[\mathscr M :
\zp[\Gamma_n]T_\Delta(p \zeta_p + \widetilde{\mathscr T}_n)]$.
Note that $T_\Delta p \zeta_p \not\in \zp[\Gamma_n]T_\Delta(p
\zeta_p + \widetilde{\mathscr T}_n)$. From $$p^n T_\Delta p
\zeta_p=p^n e_0(T_\Delta(p \zeta_p + \widetilde{\mathscr T}_n))
\in \zp[\Gamma_n]T_\Delta(p \zeta_p + \widetilde{\mathscr T}_n)$$
we deduce that the index is less or equal to $p^n$. Let $p^h$ the
order of $T_\Delta p \zeta_p$. Then $p^h T_\Delta p \zeta_p=u_n
T_\Delta(p \zeta_p + \widetilde{\mathscr T}_n)$ with $u_n \in
\zp[\Gamma_n]$. This implies that $p^h e_0 \in \zp[\Gamma_n ]$ and
$h\geq n$. Thus we have $[\mathscr M : \zp[\Gamma_n]T_\Delta(p
\zeta_p + \widetilde{\mathscr T}_n)]=p^n$ which proves the
following result.

\begin{theorem}The $\zp[\Gamma_n]$-module $T_\Delta \log_p U_n$ is generated by
$T_\Delta \widetilde{\mathscr T}_n$ and $p$.
\end{theorem}

\begin{Cor} The $\zp[\Gamma_n]$-module $T_\Delta U_n$ is generated by
$T_\Delta \overline E_n$ and $(1+p)$.

In particular, let $U_n'=\{u\in U_n \,|\ N_{K_n/\mathbb
Q_p}(u)=1\}$. Then $T_\Delta \log_p U_n'=T_\Delta \log_p\overline
E_n = T_\Delta \log_p \overline C_n$.
\end{Cor}

%
%
%
%
%
%
%

\section{A result \`a la Iwasawa in the odd part}
We recall that $K_n=\mathbb Q(\zeta_{p^{n+1}})$. Let $\pi_n =
\zeta_{p^{n+1}}-1$. In this section $\epsilon_n(\theta)$ is the
element of $\limproj \zp[\Gamma_n]$ corresponding to the power
series $f(\frac{1}{1+T}-1, \omega \theta)$ where $\theta$ is an
odd character of $\Delta$. We assume that $\theta \neq
\omega^{-1}$.
\subsection{The main result}

\begin{lemma}
Let $\chi$ be an odd character in $\widehat{\Delta}$, with $\chi \not=
\omega^{-1}$. We have $e_\chi \frac{1}{\pi_n} \in
\mathbb Z_p[\zeta_{p^{n+1}}]$.
\end{lemma}
\begin{proof} It is enough to show that $\pi_n e_\chi \frac{1}{\pi_n} \equiv 0
\bmod{\pi_n}$. A straightforward calculation gives
\begin{eqnarray*} (p-1) \pi_n e_\chi \frac{1}{\pi_n} &=& \sum_{\sigma
\in
\Delta} \frac{\overline{\chi}(\sigma)
\pi_n}{\sigma(\pi_n)} \\ & \equiv & \sum_{\sigma  \in \Delta}
\omega^{-1}
(\sigma) \overline{\chi}(\sigma) \bmod{\pi_n} \\
& \equiv & 0 \bmod{\pi_n}
\end{eqnarray*}
because $\chi \not= \omega^{-1}$, which proves the lemma.
\end{proof}

Recall that $$\epsilon_n(\chi)=\frac{1}{p^{n+1}}
\sum_{\begin{matrix}1 \leq a \leq p^{n+1} \\p\not\| a\end
{matrix}}a \theta\omega^{-1}(a) \gamma_n(a)$$ where $\gamma_n(a)$
is as defined in the introduction.
\begin{lemma} Let $T_n \in K_n$ be the sum $T_n= \sum_
{d=0}^n \zeta_{p^{d+1}}$. We have $$e_\chi \frac{1}{\pi_n}=\epsilon_n(\chi)
e_\chi T_n.$$
\end{lemma}
\begin{proof}Let $f(X)=\frac{X^{p^{n+1}}-1}{X-1}$. It is easily checked that
$$ Xf'(X)=\sum_{k=1}^{p^{n+1}-1} kX^k \quad \textrm{and} \quad
(X-1)Xf'(X)+Xf(X) =p^{n+1}X^{p^{n+1}}.$$ Letting $X =
\zeta_{p^{n+1}}$ we obtain
\begin{eqnarray*}
\frac{1}{\pi_n}&=&\frac{1}{p^{n+1}}\sum_{k=1}^{p^{n+1}-1} k
\zeta_{p^{n+1}}^k \\ & = & \frac{1}{p^{n+1}}
\sum_{\begin{matrix} \scriptstyle a=1 \\ \scriptstyle p\not \| a\end{matrix}}^
{p^{n+1}} a \zeta_{p^{n+1}}^a + \frac{1}{p^n}
\sum_{k=1}^{p^n-1}k \zeta_{p^n}^k.
\end{eqnarray*}
By induction we get $e_\chi \frac{1}{\pi_n}=\sum_{k=1}^n
\epsilon_k(\chi) e_\chi(\zeta_{p^{k+1}})$. For $m \leq n$ the
restriction of $\epsilon_n(\chi)$ to $K_m$ is $\epsilon_m(\chi)$.
The lemma follows.
\end{proof}

We now deal with the case where $\chi=\omega$. Consider the power
series $g(T)=(1-\frac{1+p}{1+T})f(T,1)$. We have
$$g(\frac{1}{1+T}-1)=(1-(1+p)(1+T))f(\frac{1}{1+T}-1,1).$$ Let
$(\epsilon_n)_{n \geq0}$ correspond to $g(\frac{1}{1+T}-1)$. The
same calculation as in lemma 2.3 shows that for all $n \geq 0$
$$(1-(1+p)\gamma_0) e_{\omega^{-1}} \frac{1}{\pi_n}=\epsilon_n
e_{\omega^{-1}} T_n.$$

Computing the sum $\sum_{\chi odd} \epsilon_n e_\chi
(T_n-\overline T_n)$ and using theorem 1.10 we easily obtain the
following result.
\begin{theorem}
There exists $\nu_\infty=(\nu_n)_{n \geq0} \in \limproj U_n^-=U_\infty^-$ such that
each $\nu_n \in U_n^-$ is unique modulo
$\mu_{p^{n+1}}$ and  $$ \log_p \nu_n = \sum_{i=0}^n p^{i-n}
e_d (\frac{1}{\pi_n}-\frac{1}{\overline \pi_n}-2(1+p)\gamma_0
e_{\omega^{-1}} \frac{1}{\pi_n}).$$ 
For all odd characters $\chi \in
\widehat{\Delta}\setminus\{\omega\}$ we have $$e_\chi
U_\infty^-/\Lambda e_\chi \nu_\infty \simeq \Lambda/ \widetilde f
(\frac{1}{1+T}-1, \omega\chi) \Lambda$$ where $ \widetilde f
(\frac{1}{1+T}-1, \omega\chi) =  f (\frac{1}{1+T}-1, \omega\chi)$
when $\chi \neq \omega^{-1}$ and $ \widetilde f (\frac{1}{1+T}-1,
\omega\chi)= (1-(1+p)(1+T)) f (\frac{1}{1+T}-1,1)$ when $\chi =
\omega^{-1}$.\\ For $\chi=\omega$ we have $$e_\chi \log_p
U_\infty^-/\Lambda e_\chi \log_p \nu_\infty \simeq \Lambda/ f
(\frac{1}{1+T}-1, \omega^2) \Lambda.$$
\end{theorem}


\subsection{A result \`a la Stickelberger}
The aim of this section is to obtain a global result from theorem
2.3 which is of local nature. In order to achieve this, set
$$\epsilon_n=\frac{1}{p^{n+1}} \sum_{\begin{matrix} 1 \leq a \leq
p^{n+1}\\ p \not \|a
\end{matrix}} a \, \sigma_a\in \mathbb Q[G_n]$$ where
$G_n=\Delta \times \Gamma_n$. This element looks like the
Stickelberger one and our computation are inspired by this
analogy. Notice that the restriction of $\epsilon_{n+1}$ to
 $K_n$ is not $\epsilon_n$ but
$\epsilon_n+(p-1)/2 N_n$ where $N_n$ is the norm element of the
group algebra $\mathbb Z[G_n]$. We should thus consider the
element $(j-1)\epsilon_n$ where $j$ is the complex conjugation
which is compatible with the canonical morphisms $\mathbb
Z[G_{n+1}]~\rightarrow~\mathbb Z[G_n]$. By lemma 2.2 we have
$\frac{1}{\pi_n}=\sum_{k=0}^n \epsilon_k(\zeta_{p^{n+1}})$. We
deduce that $(j-1)\frac{1}{\pi_n}=(j-1) \epsilon_n T_n$.

\paragraph{} Let $\theta_n=\frac{1}{\pi_n}-\frac{1}{\overline{\pi}_n}$.
Define the ideal $I$ of $\mathbb Z[G_n]$ by $I=\mathbb
Z[G_n](\sigma_c-c^*)$ where $c^*$ is the inverse of $c$ modulo
$p^{n+1}$. However contrary to the standard case (the one of the
Stickelberger element), $I$ is not the order associated to
$\frac{1}{\pi_n}$. We also define the ideal $\mathscr{I}=\mathbb
Z[G_n]\epsilon_n \cap \mathbb Z[G_n]$.

\paragraph{} We need two more ideals.
$$\left\{ \begin{matrix}\mathcal{E}_n & = & \mathbb Z[G_n] T_n
\mbox{ and}\\ \mathcal{C}_n & = & I\theta_n \end{matrix} \right.$$
Notice that $\mathcal C_n = \mathcal C_n^-$. The aim of this
section is to show the following result.
\begin{theorem}We have $[\,\mathcal{E}_n^-:\mathcal C_n\,] =
2^{\frac{|G_n|}{2}-1}\cdot h_{p^{n+1}}^-$ where $h_{p^{n+1}}^-$ is
the quotient of the class number of $K_n$ by the class number of
the maximal real subfield of $K_n$.
\end{theorem}
\begin{proof} We need several lemmas.
\begin{lemma} We have $\mathscr{I}=I \epsilon_n$.
\end{lemma}
\begin{proof} We want to show that for all $\beta \in \mathbb Z[G_n]$,
 $\beta \epsilon_n \in \mathbb Z[G_n]$ is equivalent to
 $\beta \in I$.
Let $\beta \in I$. We can assume that $\beta=(\sigma_c-c^*)$. Then
\begin{eqnarray*} (\sigma_c-c^*) &=&\sum_a \{\frac{a}{p^{n+1}}\}
\sigma_{ac}-c^*\{\frac{a}{p^{n+1}}\} \sigma_a \\ &=& \sum_a
(\{\frac{ac^*}{p^{n+1}}\} -c^* \{\frac{a}{p^{n+1}}\})\sigma_a
.\end {eqnarray*} As $cc^*\equiv 1 \bmod{p^{n+1}}$ we obtain
$(\sigma_c-c^*)\epsilon_n \in \mathbb Z[G_n]$.Conversely, notice
first that $p^{n+1}=p^ {n+1}-1+\sigma_{1+p^{n+1}} \in I$. Let
$\beta= \sum_a x_a \, \sigma_a$ with $x_a \in \mathbb Z$ and
assume that $\beta \epsilon_n \in \mathbb Z[G_n]$. Then
\begin{eqnarray*} \beta \cdot \epsilon_n &=& \sum_a\sum_c x_a
\{\frac{c}{p^ {n+1}}\}\sigma_{ac} \\ &=& \sum_b \sum_a x_a
\{\frac{a^*b}{p^{n+1}}\} \sigma_b.
\end{eqnarray*}
When $b=1$ our assumption implies that $p^{n+1}|\sum_a x_a a^*$ so
that $\sum_a x_a a^* \in I$. Finally we get $$ \beta =\sum x_a a^*
= \sum x_a (\sigma_a-a^*)+\sum x_a a^* \in I.$$

\end{proof}
 Let us return to the proof of theorem 2.4. On the one hand, by
 Leopoldt's theorem we know that $\mathcal E_n$ is a free $\mathbb Z[G_n]$-module
of rank one. We thus have $\mathcal E_n^- = \mathbb Z[G_n]^-
~\cdot~T_n$. A straightforward calculation shows that $\mathbb
 Z[G_ n]^-=(j-1)\mathbb Z[G_n]$. On the other hand, we have  $\theta-n=(j-1)
\epsilon_n T_n$ (see beginning os this section). Then we get
 $$\mathcal C_n= (j-1) \mathscr{I} \epsilon_n T_n=(j-1) \mathscr(I) T_n$$
In order to complete the proof of the theorem we have to express
$\mathscr I^-$ according to $(j-1) \mathscr I$ and to compute the
index of $\mathscr I ^-$ in $\mathbb Z[G_n]^-$. By definition we
have $$\mathscr I^- = \mathscr I \cap \mathbb Z[G_n]^- = \mathbb
Z[G_n] \cdot \epsilon_n \cap \mathbb Z[G_n]^-.$$ The following
result resembles a theorem by Iwasawa (see \cite{WAS} theorem
6.19).
\begin{prop}We have
$[\,\mathbb Z[G_n]^- : \mathscr I^- \,] = h_{p^{n+1}}^-$.
\end{prop}
\begin{proof}The proof is the same as the one of Iwasawa's theorem.
Let us recall the main steps : first complete and then worke at
each prime. We define the ideal $\mathscr{I}_q = \mathbb Z_q [G_n]
\mathscr{I}$ for a prime $q$. We have the following results.
\begin{lemma}
\begin{enumerate}
\item $\mathbb Z_q [G_n]^-=(1-j) \mathbb Z_q [G_n]$
\item $ \mathscr{I}_q= \mathbb Z_q [G_n] \cdot \epsilon_n \cap \mathbb Z_q
[G_n]$
\item $ \mathscr{I}_q^-=\mathscr{I}_q \cap \mathbb Z_q [G_n]=\mathbb Z_q [G_n]
\cdot \epsilon_n \cap \mathbb Z_q [G_n]^-$
\item $\mathscr{I}_q^- = \mathscr I^- \cdot \mathbb Z_q $
\item When $p \not= q$, $\mathscr I_q = \mathbb Z_q [G_n] \epsilon_n$
\end{enumerate}
\end{lemma}
The proof runs just as in \cite{WAS}, lemma 6.20. By lemma 2.7 we
have an isomorphism $\mathbb Z_q [G_n]^-/\mathscr I_q^- \simeq
(\mathbb Z [G_n]^-/\mathscr I^-) \otimes \mathbb Z_q$. It is
enough to show the following result.
\begin{prop} The index $[\, \mathbb Z_q [G_n]^- : \mathscr I_q^-
\,]$ is the q-part of  $h_{p^{n+1}}^-$ for all primes $q$.
\end{prop}
\begin{proof}
Assume to begin with that $q\not= 2,p$. Then $(1 \pm j)/2 \in
\mathbb Z_q[G_n]$ so we can separate the plus-part and the
minus-part. We obtain $$\mathscr I_q^-=\frac{1-j}{2} \mathscr I =
\mathbb Z_q[G_n]^- \cdot \epsilon_n$$ We have to calculate the
index $[\, \mathbb Z_q[G_n]^-: \mathbb Z_q[G_n]^- \cdot \epsilon_n
\, ]$ which equals to the $q$-part of the determinant of the map
$$
\begin{array}{rcl}\varphi : \mathbb Z_q[G_n]^- & \rightarrow &
\mathbb Z_q[G_n]^- \\ x &\mapsto& x\, \epsilon_n.
\end{array}$$ Compute this determinant in $\overline{\mathbb
Q}_q[G_n]^- = \oplus_{\chi \textrm{ odd}} \,e_\chi
\overline{\mathbb Q}_q[G_n]$. However we have $e_\chi
\sigma=\chi(\sigma) e_\chi$ for all $\sigma \in G$ from which we
deduce that $$ e_\chi \epsilon_n = B_{1, \chi} e_\chi.$$ Then we
have
\begin{eqnarray*} [\, \mathbb Z_q[G_n]^- : \mathscr I_q^-\,] &
=& q-\mbox{ part of det($\varphi$)}=q-\mbox{part of } \prod_{\chi
\textrm{ odd}} B_{1,\chi}\\ &=& q-\mbox{part of }h_{p^{n+1}}^-.
\end{eqnarray*}
\\ \indent Let us now deal with the case $p=2$. The trick is to modify $\epsilon_n$
 and to define
$\tilde{\epsilon}_n=\epsilon_n-1/2 N$ where $N$ is the norm
element of the group algebra $\mathbb Z[G_n]$. We easily check
that $\frac{1-j}{2} \widetilde{\epsilon}_n= \tilde{\epsilon}_n$ so
that $\tilde{\epsilon}_n \in \mathbb Q_2[G_n]^-$.
\begin{lemma}We have \begin{enumerate}
\item $\mathscr I_2^- \subseteq \mathbb Z_2[G_n]
\tilde{\epsilon}_n $
\item $[\, \mathbb Z_2[G_n] \tilde{\epsilon}_n :\mathscr{I}_2^-
\,]=2$
\end{enumerate}
\end{lemma}
\begin{proof}
The first part is obvious. For the second one, notice that if
 $x\in \mathbb Z_2[G_n]$ then either $x \,
\tilde{\epsilon}_n \in \mathbb Z_2[G_n]$ or $x \,
\tilde{\epsilon}_n-\tilde{\epsilon}_n \in \mathbb Z_2[G_n]$, and
that $ \mathbb Z_2[G_n] \tilde{\epsilon}_n \cap \mathbb Z_2[G_n]=
\mathscr{I}_2$.
\end{proof}
The end of the proof for $p=2$ runs the same as previously except
that the map $\varphi$ is now the multiplication by
$\tilde{\epsilon}_n$. Moreover we have $e_\chi \tilde{\epsilon}_n=
e_\chi \epsilon_n$ when $\chi$ is odd. Then we get $$[\,\mathbb
Z_2[G_n]^-:\mathbb Z_2[G_n]^- \tilde{\epsilon}_n\,]= 2-\mbox{part
of } \det \varphi = 2^{(1/2)|G_n|-1} ( 2-\mbox{part of }
h_{p^{n+1}}^-).$$ Note that $$ \mathbb Z_2[G_n]^-
\tilde{\epsilon}_n=(1-j)\mathbb Z_2[G_n] \tilde{\epsilon}_n =
\mathbb Z_2[G_n] (2\tilde{\epsilon}_n)=2 \mathbb Z_2[G_n]
\tilde{\epsilon}_n.$$ It follows that $[\, \mathbb Z_2[G_n]
\tilde{\epsilon}_n : \mathbb Z_2^-[G_n]
\tilde{\epsilon}_n\,]=2^{\mathbb Z_2-\textrm{rank of } \mathbb
Z_2[G_n]^-}= 2^{(1/2)|G_n|}$. Together with lemma 2.6 this yields
the required result $$[\, \mathbb Z_2[G_n]^-:\mathscr
I_2^-\,]=2-\mbox{part of }h_{p^{n+1}}^-.$$

 \indent It remains to deal with the case where $q=p$. Let us consider again the element
 $\tilde{\epsilon}_n=\epsilon_n -
\frac1 2 N$. Notice that as in \cite{WAS} we have the equivalence
$x \tilde{\epsilon}_n \in \mathbb Z_p[G_n]^- \iff x \epsilon_n \in
\mathbb Z_p[G_n]$ for $x \in \mathbb Z_p[G_n]$. The equality $[\,
\mathbb Z_p[G_n]\tilde{\epsilon}_n : \mathbb Z_p[G_n]
\tilde{\epsilon}_n \cap \mathbb Z_p[G_n]^-\,]=p^n$ follows.
However we also have $$\mathbb Z_p[G_n] \tilde{\epsilon}_n \cap
\mathbb Z_p[G_n]^- = \mathscr I_p^- \quad \textrm{and} \quad
\mathbb Z_p[G_n] \tilde{\epsilon}_n = \mathbb Z_p[G_n]^-\epsilon_n
,$$ thus $$[ \,\mathbb Z_p[G_n]^-\epsilon_n :\mathscr I_p^-
\,]=p^n.$$ We define the map $$\varphi : \mathbb Z_p[G_n]^-
\rightarrow \mathbb Z_p[G_n]^-, \quad x \mapsto p^n \epsilon_n
x.$$ Then $[\,\mathbb Z_p[G_n]^-: p^n \mathbb
Z_p[G_n]^-\epsilon_n\,]=
p^{(n/2)|G_n|}(\frac{1}{p^n})(p-\mbox{part of }h_{p^{n+1}}^-)$. We
obtain $$[\,\mathbb Z_p[G_n]^-:\mathscr I_p^-\,]=p-\mbox{ part of
}h_{p^{n+1}}^-$$ which concludes the proof of proposition 2.8.
\end{proof}
Proposition 2.6 is now proved.
\end{proof}

\begin{lemma} The index $[\,\mathscr I ^-: (1-j)\mathscr I\,]$
equals $[\mathscr I_2^- :(1-j)\mathscr I_2] =
2^{\frac{|G_n|}{2}-1}.$
\end{lemma}
\begin{proof}Localize at each prime by applying part 4 of lemma
2.7. Then
\begin{eqnarray*} [\,\mathscr I ^-: (1-j)\mathscr I\,] & =&
\prod_p [\,\mathscr I_p ^-: (1-j)\mathscr I_p \,]\\ &=&[\,\mathscr
I_2 ^-: (1-j)\mathscr I_2\,] \cdot [\,\mathscr I_p ^-:
(1-j)\mathscr I_p\,] \cdot \prod_{q\not= 2,p} [\,\mathscr I_q ^-:
(1-j)\mathscr I_q\,].
\end{eqnarray*}
When $q\neq 2,p$ lemma 2.7 implies $[\,\mathscr I_q ^-:
(1-j)\mathscr I_q\,]=1$. Moreover $\frac{1-j}{2} \in \mathbb
Z_p[G_n]$ so $\mathscr I_p^-=\frac{(1-j)}{2} \mathscr{I}_p$ and
$[\,\mathscr I_p ^-: (1-j)\mathscr I_p\,]=1$. It remains to
calculate the index at $q=2$. Notice that
$(1-j)\mathscr{I}_2=(1-j)\mathbb Z_2[G_n] \epsilon_n=2 \mathbb
Z[G_n] \tilde{\epsilon}$. We have $[\, \mathbb
Z_2[G_n]\tilde{\epsilon}_n:\mathscr{I}_2^-\,]=2$. Furthermore we
have $$[\, \mathbb Z_2[G_n]\tilde{\epsilon}_n: 2 \mathbb
Z_2[G_n]\tilde{\epsilon}_n\,]=2^{\mathbb Z_2-\mbox{rank of
}\mathbb Z_2[G_n]\tilde{\epsilon}_n}=2^{(1/2)|G_n|}$$ which
concludes the proof of the lemma and of the theorem.
\end{proof}
\end{proof}

\end{document}